\documentclass[a4paper,10pt]{article}
\usepackage[latin1]{inputenc}
\usepackage{amsthm}
\usepackage{amsmath,amssymb,amsfonts}
\usepackage{graphicx}
\usepackage{epic}
\usepackage[french,english]{babel}
\usepackage[mathscr]{euscript}
\usepackage{enumerate}
\usepackage{xspace}

\input xy
\xyoption{all}
\usepackage[dvips,all]{xy}
\input xypic
\xyoption{v2}
\usepackage[pdftex]{hyperref}
\usepackage{url}

\makeatletter
\long\def\nnfoottext#1{\insert\footins{\footnotesize
    \interlinepenalty\interfootnotelinepenalty
    \splittopskip\footnotesep
    \splitmaxdepth \dp\strutbox \floatingpenalty \@MM
    \hsize\columnwidth \@parboxrestore
   \edef\@thefnmark{}
   \edef\@currentlabel{}\@makefntext
    {\rule{\z@}{\footnotesep}\ignorespaces
      #1\strut}}}
\makeatother

\newtheorem{thm}{Theorem}

\newtheorem{lm}[thm]{Lemma}

\theoremstyle{definition}

\theoremstyle{remark}

\theoremstyle{plain}

\newcommand{\N}{\mathbb{N}}
\newcommand{\K}{\mathbb{\Bbbk}}

\newcommand{\g}{\mathfrak{g}}
\newcommand{\h}{\mathfrak{h}}

\newcommand{\lf}{\mathfrak{l}}

\newcommand{\Od}{\mathcal{O}}
\newcommand{\NN}{\mathcal{N}}

\newcommand{\sln}{\mathfrak{sl}}
\newcommand{\so}{\mathfrak{so}}

\newcommand{\gl}{\mathfrak{gl}}

\newcommand{\Aut}{\mathrm{Aut}}

\title{The closure of a sheet is not always a union of sheets, a~short note}
\author{\sc Micha\"el Bulois
\thanks{{\url{michael.bulois@univ-angers.fr}}}}
\begin{document}
\selectlanguage{english}
\date{}
\maketitle
\nnfoottext{Universit\'e d'Angers,  Laboratoire de Math\'ematiques d'Angers, UMR 6093, 2 boulevard Lavoisier
49045 Angers CEDEX 1, France}
\begin{abstract}
In this note we answer to a frequently asked question. If $G$ is an algebraic group acting on a variety $V$, a $G$-sheet of $V$ is an irreducible component of $V^{(m)}$, the set of elements of $V$ whose $G$-orbit has dimension $m$. We focus on the case of the adjoint action of a semisimple group on its Lie algebra. We give two families of examples of sheets whose closure is not a union of sheets in this setting.
\end{abstract}

Let $\g$ be a semisimple Lie algebra defined over an algebraically closed field $\K$ of characteristic zero.
Let $G$ be the adjoint group of $\g$. For any integer $m$, one defines $$\g^{(m)}=\{x\in\g\mid \dim G.x=m\}.$$
A $G$\emph{-sheet} (or simply \emph{sheet}) is an irreducible component of $\g^{(m)}$ for some $m\in\N$. We refer to \cite[\S39]{TY} for elementary properties of sheets. The most important one is that each sheet contains a unique nilpotent orbit.

There exists a well known subdivision of sheets which forms a stratification. The objects considered in this subdivision are Jordan classes and generalize the classical Jordan's block decomposition in $\gl_n$. Those classes and their closures are widely studied in \cite{Bo1} (cf. also \cite[\S39]{TY} for a more elementary viewpoint).
Since sheets are locally closed, a natural question is then the following. 
$$\mbox{If $S$ is a sheet, is $\overline{S}$ is a union of sheets?}$$
The answer is negative in general. We give two families of counterexamples below.

\begin{enumerate}
\item  A nilpotent orbit $\Od$ of $\g$ is said to be rigid if it is a sheet of $\g$. Rigid orbits are key objects in the description of sheets given in \cite{Bo1}. They were classified for the first time in \cite[\S II.7\&II.10]{Sp}.
The closure ordering of nilpotent orbits (or \emph{Hasse diagram}) can be found in \cite[\S II.8\&IV.2]{Sp}. In the classical cases, a more recent reference for these lists is \cite{CM}.
One easily checks from these classifications that there may exists some rigid nilpotent orbit $\Od_{1}$ that contains a non-rigid nilpotent orbit $\Od_{2}$ in its closure.
Then, we set $S=\Od_{1}$ and we get $\Od_{2}\subset\overline{S}\subset\NN(\g)$ where $\NN(\g)$ is the set of nilpotent elements of $\g$. 
Since $\Od_{2}$ is not rigid, the sheets containing $\Od_{2}$ are not wholly included in $\NN(\g)$. 
Therefore, the closure of $S$ is not a union of sheets.

Here are some examples of such nilpotent orbits.
In the classical  cases, we embed $\g$ in $\gl_n$ in the natural way. Then, we can assign to each nilpotent orbit $\Od$, a partition of $n$, denoted by $\Gamma(\Od)$. This partition defines the orbit $\Od$, sometimes up to an element of $\Aut(\g)$.
In the case $\g=\so_{8}$ (type D$_{4}$), there is exactly one rigid orbit $\Od_{1}$, such that $\Gamma(\Od_1)=[3,2^2,1]$. It contains in its closure the non-rigid orbit $\Od_{2}$ such that $\Gamma(\Od_2)=[3,1^5]$ (cf.~\cite[Table2, p.15]{Mo}). 
Very similar examples can be found in types C and B. 

In the exceptional cases, we denote nilpotent orbits by their Bala-Carter symbol as in \cite{Sp}. 
Let us give some examples of the above described phenomenon. 
\begin{itemize}
\item in type $E_6$ ($\Od_1=3A_1$ and $\Od_{2}=2A_1$), 
\item in type $E_7$ ($\Od_1=A_2+2A_1$ and $\Od_{2}=A_2+A_1$), 
\item in type $E_8$ ($\Od_1=A_2+A_1$ and $\Od_{2}=A_2$) 
\item and in type $F_4$ ($\Od_1=A_2+A_1$ and $A_2$).   
\end{itemize}

\item In the case $\g=\sln_{n}$ of type $A$, there is only one rigid nilpotent orbit, the null one. Hence the phenomenon depicted in 1 can not arise in this case. Let $S$ be a sheet and let $\lambda_S=(\lambda_{1}\geqslant\dots\geqslant \lambda_{k(\lambda_S)})$ be the partition of $n$ associated to the nilpotent orbit $\Od_S$ of $S$. 
As a consequence of the theory of induction of orbits, cf. \cite{Bo1}, we have \begin{equation}\label{eq}\overline{S}=\overline{G.\h_{S}}^{reg}\end{equation} where $\h_{S}$ is the centre of a Levi subalgebra $\lf_{S}$. The size of the blocks of $\lf_{S}$ yield a partition of $n$, which we denote by $\tilde{\lambda}_S=(\tilde{\lambda}_{1}\geqslant\dots\geqslant \tilde{\lambda}_{p(\lambda_S)})$. In fact $\tilde{\lambda}$  is the dual partition of $\lambda$, i.e. $\tilde{\lambda_{i}}=\#\{j\mid \lambda_{j}\geqslant i\}$ (see, e.g.,  \cite[\S2.2]{Kr}). In particular, the map sending a sheet $S$ to its nilpotent orbit $\Od_S$ is a bijection.

An easy consequence of \eqref{eq} is the following (see \cite[Satz 1.4]{Kr}). Given any two sheets $S$ and $S'$ of $\g$, 
we have $S\subset \overline{S'}$ if and only if $\h_{S}$ is $G$-conjugate to a subspace of $\h_{S'}$  or, equivalently, $\lf_{S'}$ is conjugate to a subspace of $\lf_{S}$. This can be translated in terms of partitions by defining a partial ordering on the set of partitions of $n$ as follows.
We say that $\lambda\preceq\lambda'$ if there exists a partition $(J_i)_{i\in [\![1,p(\lambda)]\!]}$ of $[\![1,p(\lambda')]\!]$ such that $\tilde{\lambda}_i=\sum_{j\in J_i} \tilde{\lambda}'_j$. Hence, we have the following characterization.
\begin{lm}
$S\subset \overline{S'}$ if and only if $\lambda_S\preceq\lambda_{S'}$.
\end{lm}
One sees that this criterion  is strictly stronger than the characterization of inclusion relations of closures of nilpotent orbits (see, e.g., \cite[\S6.2]{CM}).
More precisely, one easily finds two sheets $S$ and $S'$ such that $\Od_{S}\subset\overline{\Od_{S'}}$ while $\lambda_{S}\npreceq\lambda_{S'}$.
Then, $\Od_S\subset \overline{S'}$, $S$ is the only sheet containing $\Od_S$ and $S\not\subset \overline{S'}$.
For instance, take $\lambda_{S'}=[3,2]$, $\lambda_{S}= [3,1,1]$. Their respective dual partitions being $[2,2,1]$ and $[3,1,1]$, we have $\lambda_{S}\npreceq\lambda_{S'}$.
\end{enumerate}


\end{document}